\documentclass[a4paper,UKenglish,cleveref, autoref,nosubfigcap, nolineno]{article}
\usepackage[a4paper,
            bindingoffset=0.2in,
            left=1.0in,
            right=1.0in,
            top=1.0in,
            bottom=1.0in,
            footskip=.25in]{geometry}
\usepackage[usenames, dvipsnames]{xcolor}
\usepackage{tikz-cd}
\usepackage{tikz}
\usepackage{adjustbox}
\usepackage[normalem]{ulem}
\usepackage{multirow}
\usepackage{amsmath}
\usepackage{amsthm,thmtools}
\usepackage{hyperref}

\theoremstyle{plain} 
\newtheorem{theorem}{Theorem}[section]

\newtheorem{prop}{Proposition}[section]

\theoremstyle{definition}
\newtheorem{defn}{Definition}[section]

\theoremstyle{definition}

\theoremstyle{remark}

\usepackage[linesnumbered,algoruled,boxed,lined,noend]{algorithm2e}

\def\corcommstyle#1{\bf\small\tt\textcolor{Mahogany}{#1}}
\def\corcommredstyle#1{\bf\small\tt\textcolor{Mahogany}{#1}}

\def\corrl #1<<#2||#3>>{
\if\visiblecomments y
  \begin{quote} {\corcommstyle{ $<<$COMMENT$>>$ #1\marginpar{!!}\\$<<$OLD$<<$}} \end{quote}
  #2
  \begin{quote} {\corcommstyle{ ==NEW== }} \end{quote}
  #3
  \begin{quote} {\corcommstyle{ $>>$END$>>$ }} \end{quote}
 \else
  #3
 \fi
}

\long\def\longcorrl #1<<#2||#3>>{
\if\visiblecomments y
  \begin{quote} {\corcommstyle{ $<<$COMMENT$>>$ #1\marginpar{!!}\\$<<$OLD$<<$}} \end{quote}
  #2
  \begin{quote} {\corcommstyle{ ==NEW== }} \end{quote}
  #3
  \begin{quote} {\corcommstyle{ $>>$END$>>$ }} \end{quote}
 \else
  #3
 \fi
}

\def\corrq #1<<#2>>{
\if\visiblecomments y
  \begin{quote} {\corcommstyle{ $<<$COMMENT$>>$ #1\marginpar{!!}\\$<<$BEG$<<$}} \end{quote}
  #2
  \begin{quote} {\corcommstyle{ $>>$END$>>$ }} \end{quote}
 \else
  #2
 \fi
}

\long\def\longcorrq #1<<#2>>{
\if\visiblecomments y
  \begin{quote} {\corcommstyle{ $<<$COMMENT$>>$ #1\marginpar{!!}\\$<<$BEG$<<$}} \end{quote}
  #2
  \begin{quote} {\corcommstyle{ $>>$END$>>$ }} \end{quote}
 \else
  #2
 \fi
}

\def\corrc #1<<>>{
\if\visiblecomments y
  \begin{quote} {\corcommstyle{ $<<$COMMENT$>>$ #1\marginpar{!!}}} \end{quote}
\fi
}

\def\corrse #1<<>>{
\if\visiblecomments y
  \begin{quote} {\corcommstyle{ $<<$SECOND EDITION$>>$ #1\marginpar{!!}}} \end{quote}
\fi
}

\def\corre #1<<#2||#3>>{
\if\visiblecomments y
  #3\marginpar{\corcommstyle{ #1}}
 \else
  #3
 \fi
}

\long\def\longcorre #1<<#2||#3>>{
\if\visiblecomments y
  #3\marginpar{\corcommstyle{ #1}}
 \else
  #3
 \fi
}

\def\corrs #1<<#2||#3>>{
\if\visiblecomments y
  \textcolor{Mahogany}{#3}\marginpar{\corcommstyle{ #2 $\rightarrow$ #3\\ #1}}
 \else
  #3
 \fi
}

\def\corrm #1<<#2>>{
\if\visiblecomments y
  \textcolor{blue}{#2}\marginpar{\corcommstyle{#1}}
 \else
  #2
 \fi
}

\def\corrmr #1<<#2>>{
\if\visiblecomments y
  \textcolor{blue}{#2}\marginpar{\corcommredstyle{#1}}
 \else
  #2
 \fi
}

\def\corro #1<<#2||#3>>{
#2}

\def\corrn #1<<#2||#3>>{
#3}

\long\def\longcorro #1<<#2||#3>>{
#2}

\long\def\longcorrn #1<<#2||#3>>{
#3}

\long\def\underconstruction #1<<<#2>>>{
\if\visiblecomments y
  \begin{quote} {\corcommstyle{ $<<$UNDER CONSTRUCTION - BEGIN$>>$ #1\marginpar{!!}}} \end{quote}
  #2
  \begin{quote} {\corcommstyle{ $>>$UNDER CONSTRUCTION - END$>>$ }} \end{quote}
 \else
 \fi
}

\def\showcomments{
  \let\visiblecomments y
}

\def\hidecomments{
  \let\visiblecomments n
}

\let\visiblecomments y 
\usepackage{amssymb}

\def\refeq#1{\if\workingver y(\ref{#1})-[[#1]]\else(\ref{#1})\fi}
\def\refth#1{\if\workingver y\ref{#1}-[[#1]]\else\ref{#1}\fi}
\def\mylabel#1{\if\workingver y\label{#1}{\bf\ \ [[#1]]\ \ }\else\label{#1}\fi}
\def\mybibitem#1{\if\workingver y\bibitem{#1}{\bf\ \ [[#1]]\ \
}\else\bibitem{#1}\fi}

\renewcommand{\emptyset}{\varnothing}
\renewcommand{\rho}{\varrho}
\renewcommand{\phi}{\varphi}
\renewcommand{\epsilon}{\varepsilon}

\def\cM{\text{$\mathcal M$}}

\def\cV{\text{$\mathcal V$}}

\def\bE{\text{$\mathbf E$}}

\def\be{\text{$\mathbf e$}}

\newcommand{\id}{\operatorname{id}}
\newcommand{\cl}{\operatorname{cl}}

\renewcommand{\emptyset}{\varnothing}

\def\proof{{\bf Proof:\ }}

\def\begeq#1{\begin{equation}\mylabel{#1}}
\def\endeq{\end{equation}}

\def\mathobj#1{\mbox{$#1$}}

\def\II{\mathobj{\mathbb{I}}}

\def\PP{\mathobj{\mathbb{P}}}

\def\ZZ{\mathobj{\mathbb{Z}}}

\def\scalprod#1{\langle #1 \rangle}

\def\implies{\;\Rightarrow\;}

\def\setof#1{\mbox{$\{\,#1\,\}$}}

\def\0#1{\hbox{\kern25pt}$ #1 $\\}
\def\1#1{\hbox{\kern40pt}$ #1 $\\}
\def\2#1{\hbox{\kern55pt}$ #1 $\\}
\def\3#1{\hbox{\kern70pt}$ #1 $\\}

\newcounter{li}

\def\begalg#1{\begin{algo}\mylabel{#1}\normalshape:\small\baselineskip 10pt\\}
\def\endalg{\end{algo}}

\def\Figures(include=#1,cat=#2){
  \renewcommand{\textfraction}{.20}
  \renewcommand{\topfraction}{.80}
  \renewcommand{\bottomfraction}{.80}
  \renewcommand{\floatpagefraction}{.80}
  \newcount\figcount
  \figcount=0
  \let\includefigures=#1
  \def\figcat{#2}
}

\def\FigureFromFile[#1][#2](#3)#4
{
  \begin{figure}[htbp]
     \global\advance\figcount by 1
     \if\includefigures y\special{anisoscale #1.wmf, \the\hsize #2}\fi
     \vspace{#2}
     \caption{#4}
     \mylabel{#3}
   \end{figure}
}

\def\FigureFromFileTwoD[#1][#2,#3](#4)#5
{
  \begin{figure}[htbp]
     \global\advance\figcount by 1
     \if\includefigures y\special{anisoscale #1.wmf, #2 #3}\fi
     \vspace{#2}
     \caption{#5}
     \mylabel{#4}
   \end{figure}
}

\def\FigureF<#1>[#2](#3)#4
{
  \begin{figure}[htbp]
     \global\advance\figcount by 1
     \if\includefigures y\special{anisoscale \figcat/fig\number\figcount.wmf,
       \the\hsize #2}
     \fi
     \if\includefigures p
       \leavevmode
       \epsfxsize=\hsize
       \epsffile{#1}
     \fi
     \if\includefigures y
          \vspace{#2}
     \fi
     \caption{#4}
     \mylabel{#3}
   \end{figure}
}

\def\Figure[#1](#2)#3
{
  \begin{figure}[htbp]
     \global\advance\figcount by 1
     \if\includefigures y\special{anisoscale \figcat/fig\number\figcount.wmf,
       \the\hsize #1}
     \fi
     \if\includefigures p
       \leavevmode
       \epsfxsize=\hsize
       \epsffile{fig\number\figcount.eps}
     \fi
     \if\includefigures y
          \vspace{#1}
     \fi
     \caption{#3}
     \mylabel{#2}
   \end{figure}
}

\def\0{\hbox{\kern5pt}}
\def\1{\hbox{\kern20pt}}
\def\2{\hbox{\kern35pt}}
\def\3{\hbox{\kern50pt}}
\def\4{\hbox{\kern65pt}}
\def\5{\hbox{\kern80pt}}
\def\6{\hbox{\kern95pt}}

\def\kw#1{\textbf{#1}}

\def\kwif{\kw{if}\;}

\def\kwfor{\kw{for}\;}

\newcommand{\inv}{\sf Inv}

\newcommand{\gr}{\sf gr}
\newcommand{\mv}{\mathcal{V}}

\newcommand{\pgrad}[2]{{\gr}_{#1}(#2)}
\definecolor{yellow}{RGB}{255,225,55}

\newcommand{\low}{{\tt low}}

\newcommand{\leqlin}{\mathbin{\leq_{\mbox{{\scriptsize lin}}}}}
\newcommand{\PPlin}{\PP_{\mbox{{\scriptsize lin}}}}

\newcommand{\Poset}{P}
\newcommand{\md}{\mathcal{M}}

\definecolor{dark-gray}{RGB}{64,64,64}
\definecolor{medium-gray}{RGB}{114,114,114}
\definecolor{light-gray}{RGB}{190,190,190}

\begin{document}

\title{Computing Connection Matrices via Persistence-like Reductions} 

\author{Tamal K. Dey\thanks{Department of Computer Science, Purdue University, West Lafayette, Indiana, USA. \texttt{tamaldey@purdue.edu}}
\and Micha\l{} Lipi\'nski\thanks{
Dioscuri Centre in TDA, Institute of Mathematics, Polish Academy of Sciences, Warsaw, Poland\\
Division of Computational Mathematics, Faculty of Mathematics and Computer Science, Jagiellonian University, Krak\'{o}w, Poland.\texttt{michal.lipinski@uj.edu.pl}}
\and Marian Mrozek\thanks{Division of Computational Mathematics, Faculty of Mathematics and Computer Science, Jagiellonian University, Krak\'{o}w, Poland. \texttt{marian.mrozek@uj.edu.pl}}
\and Ryan Slechta\thanks{Innovation Partnerships, University of Michigan, Ann Arbor, Michigan, USA. \texttt{rslechta@umich.edu}}
}

\date{}

\maketitle
\thispagestyle{empty}

\begin{abstract}
Connection matrices are a generalization of Morse boundary operators from the classical Morse theory for gradient vector fields. Developing an efficient computational framework for connection matrices is particularly important in the context of a rapidly growing data science that requires new mathematical 
tools for discrete data.
Toward this goal, the classical theory for connection matrices has been adapted to combinatorial frameworks
that facilitate computation. We develop an efficient persistence-like algorithm to compute a connection matrix from a given combinatorial
(multi) vector field on a simplicial complex. This algorithm requires a single-pass, improving upon a known algorithm that runs an implicit recursion executing two-passes at each level.
Overall, the new algorithm is more simple, direct, and efficient than the state-of-the-art. Because of the algorithm's similarity to the persistence algorithm, one may take
advantage of various software optimizations from topological data analysis.
\end{abstract}
\section{Introduction}
Connection matrix theory, originally developed by R.\ Franzosa \cite{Fr1986,Fr1989} as a tool for proving the existence of heteroclinic
connections in dynamical systems, 
is a generalization of homological Morse theory in the setting of the Conley index theory \cite{Co78}. 
Specifically, the connection matrix is the boundary operator in a chain complex constructed from Conley indices
of Morse sets in a Morse decomposition of an isolated invariant set of a flow acting on a metrizable topological space.
The topological space replaces the smooth manifold in the classical Morse theory, the flow replaces the gradient vector
field, the Morse decomposition replaces the Morse function, and the individual Morse sets replace the critical points.
The theory has found many applications, in particular in the study of ODEs and PDEs. For some early and very recent applications, see \cite{Mo1988} and \cite{HMS2021a,MSV2021} respectively
and references therein. 

Similar to the early Conley theory, the development of applications of connection matrices has been delayed by the limitation
of manual, analytic calculations. The applicability of Conley theory rapidly changed about twenty years ago due to the development of 
modern, efficient homology algorithms. However, computation of connection matrices required more specialized algorithms. 
The first such algorithm, henceforth referred to as HMS algorithm, was presented in a recent paper by S.\ Harker, K.\ Mischaikow and K.\ Spendlove \cite{HMS2021}.
In this paper, connection matrices are studied from a computational point of view 
via a partially algebraic and partially geometric formulation
based on a simplification of the theory for field coefficients presented by J.\ Robin and D.\ Salamon in  \cite{RoSa1992}.

As pointed out by Harker, Mischaikow and Spendlove in \cite{HMS2021}, 
developing an efficient computational framework for the connection matrix theory is particularly important in the context of the rapidly growing field of data science. 
Classical mathematical models, based on differential equations, are not adequate to deal with data for which the equations cannot be
built from some first principles, like in physics, and more direct approaches are needed to address problems in dynamically changing data. 

Among such approaches is the theory developed by R. Forman, which concerns dynamical systems that are generated by combinatorial vector fields \cite{Forman1998b,Forman1998a}. Recently, various authors have shown that the Conley theory holds in the more general setting of dynamical systems generated via combinatorial multivector fields \cite{DMS2020,DMS2021,KMW2016,LKMW2022,MW2021b,Mr2017}. 

Harker, Mischaikow, and Spendlove first observed links between connection matrices and persistent homology~\cite{HMS2021}.
A connection matrix may be viewed
as a reduced matrix of a reduced filtered complex. 
In \cite{HMS2021}, the authors observed that 
a certain construction of a gradient Forman vector field $\cV$ 
consistent with the given Morse decomposition reduces the problem to constructing a connection matrix for a smaller complex
unless $\cV$ consists only of critical cells. In the latter case the connection matrix is just the matrix
of the boundary operator in the reduced complex. Hence, to compute the connection matrix they iteratively execute two-passes, one
reducing the complex (called {\sc Matching} in~\cite{HMS2021}) and another computing a
boundary operator of this reduced complex (called {\sc Gamma} in~\cite{HMS2021}) until no more reductions are possible. In a sense,
they operate \emph{recursively} starting all over again on the reduced complexes.

The partial geometric flavor of HMS algorithm is due to the fact that the Forman vector field used to reduce the problem
is constructed on the original basis acting as the phase space for the combinatorial dynamics. This helps guiding the intuition but results in the need of an implicit recursion executing two passes. 
In our algorithm we adopt this overall reduction approach while making it more efficient. In particular: 
\begin{enumerate}
    \item 
    Analogous to the matrix based persistence algorithm~\cite{EH2010}, our algorithm requires no recursive behavior as in~\cite{HMS2021}. It 
    makes a single pass over the filtered input matrix
    and performs column additions that
    implicitly execute the reductions simultaneously at both the complex and homology
    level. Then, the algorithm obtains the connection matrix by simply removing certain rows and columns. 
    A key observation that facilitates
    this single pass execution is that reductions within Morse sets and across Morse sets,
    unlike in HMS algorithm, can be combined 
    with suitable basis changes.  
    \item Like the persistence algorithm, for a multivector field on a simplicial
    complex comprising $n$ simplices,
    our single pass matrix based algorithm runs in $O(n^3)$ time whereas
    the recursive HMS algorithm runs in $O(n^4)$ time as far as we 
    can analyze (no time complexity
    analysis is given in~\cite{HMS2021}). 
    Again, the advantage ensues from embedding the homology
    computation into the reduction phase.
    \item 
    The need of implicit recursion in the HMS algorithm is closely related to the problem of the existence
    of a perfect Morse function, that is, a Morse function whose number of $n$-dimensional critical cells equals 
    the $n$th Betti number of the complex. In general, a perfect Morse function may not exist. For instance, 
    a Morse set or multivector in \cite[Figure 5]{Mr2017} has all Betti numbers zero but requires at least two critical cells 
    in a gradient Forman vector field. Also, finding a constant-factor approximation of Morse function which minimizes the number of critical cells is NP-hard \cite{LLT2003}.
    We avoid this problem 
    entirely by not forcing the algorithm to run on a fixed basis. In that sense, analogous to the persistence algorithm, 
    our algorithm is purely algebraic.
    \item Our single pass persistence-like algorithm can take advantage of
    different optimizations that have been proposed to accelerate the persistence
    algorithm in TDA, e.g.~\cite{bauer2014clear,SMV211,DW22,KS19,Wagner2012}.
\end{enumerate}



\section{Background and overview}
\subsection{Connection matrix}
Connection matrices are a generalization of Morse boundary operators in classical Morse theory for gradient vector fields. 
They represent the connections between Morse sets in the dynamics using homological algebra.
\begin{figure}[htbp]
\centerline{\includegraphics{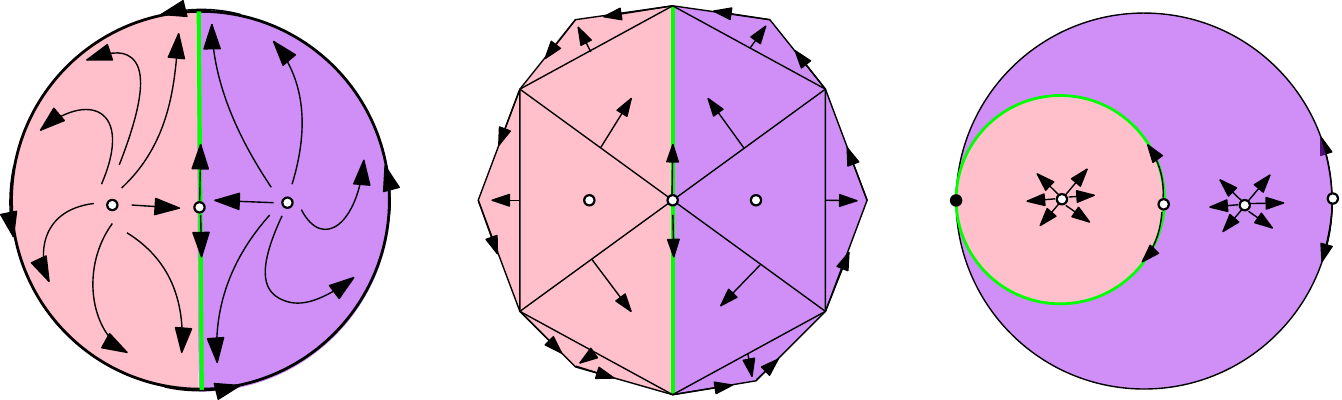}}
\caption{(left) A continuous vector field, (middle) its discretization with Forman vectors, (right) and its Conley complex represented with $5$ cells and their incidence structure.}
\label{fig:conleycomplex}
\end{figure}

\begin{figure}[htbp]
\centerline{\includegraphics{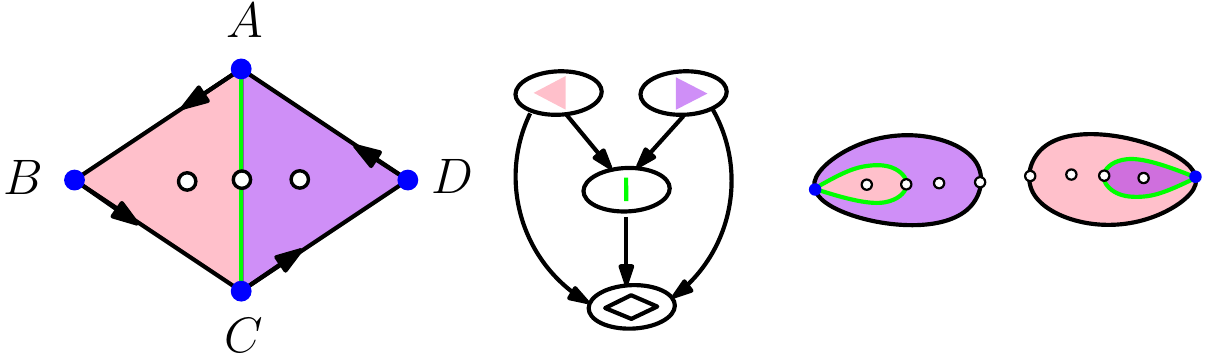}}
\caption{(left) Vector field $\mv=\{\{A,AB\},\{B,BC\},\{C,CD\},\{D,DA\},\{CA\},\{ABC\},\{CDA\}\}$, (middle) A Morse decomposition consisting of 4 Morse sets, triangles $CDA$, $ABC$, edge $CA$, and
the orbit $\{\{A,AB\},\{B,BC\},\{C,CD\},\{D,DA\}\}$, (right) Two Conley complexes; the left
one corresponds
to a connection matrix shown in Figure~\ref{fig:matrixalgo} (top) and the right one corresponds to a connection matrix shown in Figure~\ref{fig:matrixalgo} (bottom).}
\label{fig:annulus}
\end{figure}

The Conley index of a Morse set $M$ encapsulates the behavior of 
the local dynamics in terms of the homology
group of $M$ relative to its exit set. 
In particular, it differentiates between repellers, attractors and saddles.
Connection matrices, on the other hand, capture the dynamics 
at the global level by indicating the existence of connecting trajectories between Morse sets. 
In this sense, it is a useful
invariant of the complete dynamics given by the input combinatorial vector field. 
Loosely speaking,
if $H(M_p, E_p)$ denotes the Conley index for a Morse set $M_p$, that is, the homology group of $M_p$ relative to
the exit set $E_p$, then a
connection matrix represents a boundary operator
\begin{equation*}
\Delta: \bigoplus_{p\in P} H(M_p,E_p)\rightarrow \bigoplus_{p\in P} H(M_p,E_p).
\end{equation*}
Alternatively,
a connection matrix can be viewed as a simplification of the vector field producing what is called a Conley complex which is a generalization of the Morse complex studied in Morse theory.
It embodies the essential dynamics by coalescing the trivial Morse sets into larger invariant
sets with a homotopy. See e.g. Figure~\ref{fig:conleycomplex} where
a vector field with two repellers, one saddle, and one attracting orbit
is simplified with five cells, two $2$-cells representing two $2$-dimensional repellers,
one $1$-cell representing the saddle, and another $1$-cell together with
a $0$-cell representing the attracting orbit. The continuous vector field shown on the left
may be discretized with a combinatorial vector field on a simplicial complex. The Conley complex
essentially simplifies this input complex to a cell complex while preserving the
Morse sets.
Depending on how the input boundary operator at the chain level of the simplicial complex
is provided for the Morse decomposition,
the homotopy-induced simplification can be different resulting in different Conley complexes.
In Figure~\ref{fig:conleycomplex}, these different Conley complexes 
indicate the fact that the system, while preserving the Morse sets, 
may be deformed by breaking the attracting periodic orbit into an attracting stationary point and a saddle. This extra saddle may be reached from only one of the two repellers. Depending on whether this is the right or the left repeller, we get two different Conley complexes.
See, for example, the two Conley complexes in Figure~\ref{fig:annulus} (right).
Corresponding to these two Conley complexes, we have two different connection matrices
shown in Figure~\ref{fig:matrixalgo}.

\subsection{Combinatorial multivector field and Morse decomposition}
In this subsection, we introduce some of the basic definitions from combinatorial (multi)vector
field theory that are necessary for this work; see~\cite{DJKKLM19,DMS2020,LKMW2022,Mr2017} for more details.
Throughout this paper, we restrict our attention to simplicial complexes of arbitrary but finite dimension. For a simplicial complex $K$, we use $\leq$ to denote the face relation, that is, $\sigma \leq \tau$ if $\sigma$ is a face of $\tau$. We define the \emph{closure} $\cl(\sigma)$ of $\sigma$ as $\cl( \sigma ):= \{ \tau \; | \tau \leq \sigma \}$ and we extend the notion to a set of simplices $A \subseteq K$ as $\cl(A) := \cup_{\sigma \in A} \cl( \sigma )$. The set $A$ is \emph{closed} if $A = \cl(A)$. 

\begin{defn}[Multivector and multivector field]
Given a finite simplicial complex $K$, a \emph{multivector} $V\subset K$ is a subset
that is convex under face-poset relation, that is, if $\sigma,\tau\in V$ with $\sigma\leq \tau$ then every simplex $\mu$ with $\sigma\leq \mu\leq\tau$ is in $V$.
A \emph{multivector field} $\mv$ on $K$ is a partition of $K$ into multivectors.
\end{defn}

Following~\cite{LKMW2022}, a notion of dynamics on a combinatorial multivector field can be introduced. These dynamics take the form of a multivalued map $F_{\mv} \; : \; K \multimap K$, with $F_{\mv}(\sigma) = [ \sigma ]_{\mv} \cup \cl( \sigma )$ where $[ \sigma ]_{\mv}\subset K$ 
is the unique element of the partition $\cV$ containing $\sigma$. 
A finite sequence of simplices $\sigma_1, \sigma_2, \ldots, \sigma_n$ is a \emph{path} if for $i = 2, \ldots, n$, we have $\sigma_i \in F_{\mv{}}( \sigma_{i-1} )$.

\begin{defn}[Morse sets and decomposition]\label{def:ms_and_md}
Given a multivector field $\mv$ on a finite simplicial complex $K$, a subset $M\subseteq K$ 
is called
a \emph{Morse set} if 
for every path $\sigma_1, \sigma_2, \ldots, \sigma_n$ in $K$ with $\sigma_1, \sigma_n \in M$, each $\sigma_i$, $i\in \{1,\ldots,n\}$, is necessarily in $M$. A collection of Morse sets
$\md=\{M_p\,|\, p\in \Poset\}$ indexed by a poset $(\Poset, \leq)$ is called
a \emph{Morse decomposition} of $\mv$ if we have the disjoint union
$K=\sqcup_{p\in \Poset}M_p$
and $p\leq q$ if and only if there is a path
$\sigma_1,\ldots, \sigma_n$ with $\sigma_1\in M_p$ and $\sigma_n\in M_q$.
\end{defn}
An astute reader will recognize that the Morse decomposition defined here differs slightly
from the earlier definitions introduced in~\cite{DJKKLM19,DMS2021, LKMW2022}. 
Definition \ref{def:ms_and_md} liberates us from a need of introducing additional concepts
such as the invariant part $\inv\, M$ of a Morse set $M$, thereby
simplifying the presentation of the algorithm.
Nevertheless, the Morse decomposition in the sense of previous works can be easily retrieved by taking $\cM':=\{\inv\, M_p\mid p\in P\}$ (see e.g. \cite{LKMW2022} for the precise definition of the combinatorial invariance).
In particular, Definition \ref{def:ms_and_md} lets us incorporate trivial Morse sets into the structure, i.e. sets $M\in\cM$ such that $\inv\, M=\emptyset$.

\subsection{Overview of the algorithm}
Let $\mv$ be a multivector field defined on a simplicial complex $K$ with $C_*(K)$ denoting its chain spaces.
The algorithm takes the boundary operator $\partial_K: C_*(K)\rightarrow C_*(K)$
in matrix form, where the matrix is filtered according to a Morse decomposition of $\mv$.
This means that if $M_i$ and $M_j$ are two Morse sets with $i<j$, then simplices in the Morse set $M_i$ necessarily
come before those in the Morse set $M_j$ in the columns ordered from left to right.
A Morse decomposition can be obtained from $\mv$ by computing strongly connected components in a directed graph constructed as follows. Each simplex $\sigma\in K$ is represented as a node in the graph and a directed edge goes from a node $\sigma$ to a node $\tau$ if and only if $\tau \in F_\mv(\sigma)$. It is known that the strongly connected components in this graph constitute what are called the \emph{minimal Morse sets}~[Theorem 4.1,\cite{DJKKLM19}].
We do not elaborate
on this aspect any further as this is not the focus of this paper, and we assume
that the matrix $\partial_K$ is given.
Assuming coefficients in $\ZZ_2$, we perform column additions (source column added to a target column) respecting the filtration to arrive at a reduced matrix which
is analogous to the persistence algorithm; see Algorithm {\sc ConnectMat} in Section~\ref{sec:algorithm}. 
However, there are crucial differences in reducing the matrix and in outputting the result. Persistence algorithm (see e.g. books~\cite{DW22,EH2010}) takes a filtered boundary matrix representing a filtration of a simplicial complex and processes columns from left to right. Every time, the lowest entry of a column (target) `conflicts' with the lowest entry of a column (source) to the left, it adds the source column to the target column and continues such additions unless target column has its lowest entry not conflicting with the lowest entry of any column to its left or is completely zeroed out. After processing all columns, the algorithm outputs persistent simplex pairs
$(\sigma,\tau)$ if the column representing $\tau$ has the lowest entry in a row representing the simplex $\sigma$ or declares $\tau$ as unpaired if the column for $\tau$ is zeroed out.
\begin{figure}[htbp]
\centerline{\includegraphics{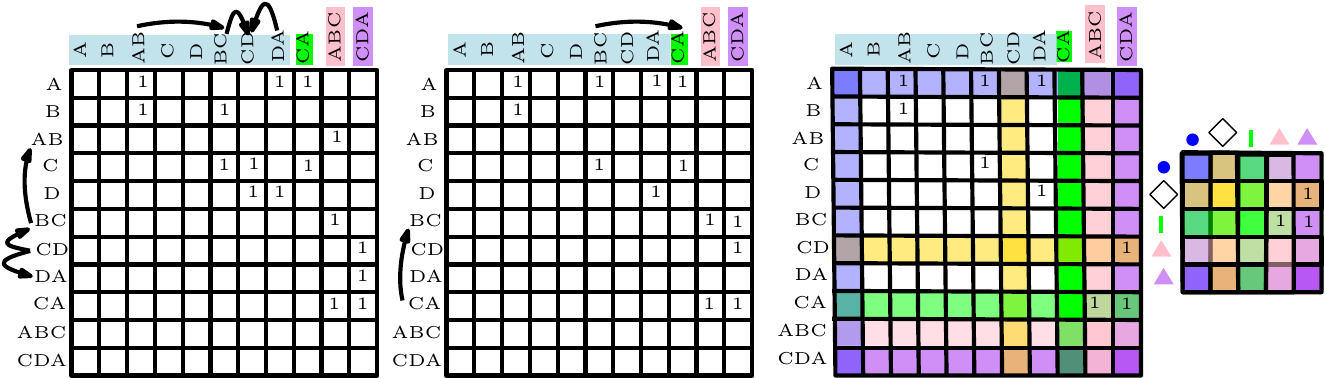}}~\\
\centerline{\includegraphics{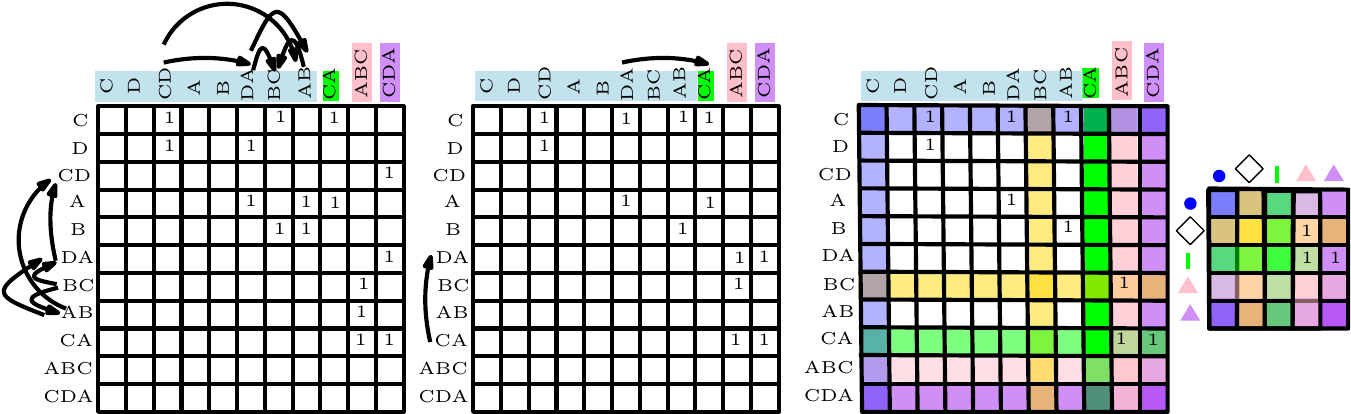}}
\caption{
Working of the algorithm on the boundary matrix for the example in Figure~\ref{fig:annulus} in two differently sorted bases; computed connection matrices are shown on right.
}
\label{fig:matrixalgo}
\end{figure}
Our algorithm {\sc ConnectMat} processes columns left to right similarly with the following differences:
 (i) column additions are performed not only to resolve conflicts at the lowest non-zero entries, but
also if the target column has a conflict for other entries with the lowest entry of the
source column, (ii) source columns are always to the left of the
target column across Morse sets, but
there is no such restrictions within a Morse set, (iii) a source column 
triggers a conflict only if it is \emph{homogeneous}, that is, 
its lowest entry comes from a simplex that is in the same Morse 
set, (iv) the algorithm also performs row additions corresponding to each column additions, (v) outputs a submatrix of the reduced matrix corresponding to columns and rows that are not homogeneous and not targetable; see discussion for the following example.

In Figure~\ref{fig:matrixalgo}, we show the boundary operator 
filtered according to the Morse decomposition
in Figure~\ref{fig:annulus}(left). It is represented by two matrices which differ only in the order of cells within the Morse sets.
We describe the working of our algorithm on the matrix shown in the top row. Also, we do not elaborate on the row operations with the implicit understanding that whenever a column $i$ is added to column $j$, the row $j$ is added to row $i$.
The column $AB$ is added to column $BC$ because of the conflict in the row $B$
even though it does not
contain the lowest $1$ in column $BC$. Similarly, column $DA$ from right is added to column $CD$ because
of a conflict in row $D$. Also, the column $BC$ from left is added to $CD$ because of a conflict
in row $C$ which renders $CD$ to be empty.
All these additions happen within the same Morse set. In the meantime, corresponding row additions affect columns $ABC$ and $CDA$. Next addition happens between
columns $BC$ and $CA$ which belong to different Morse sets, but the addition of $BC$ to $CA$ is allowed
because the lowest entry for $BC$ represents the simplex $C$ which belongs to the same Morse
set as $BC$. This column addition triggers a row addition of $CA$ to $BC$. The column $ABC$ does not trigger any conflict. The final column $CDA$ has
a conflict with the lowest entry in column $ABC$. However, $ABC$ is not added to $CDA$
because it is not homogeneous since its lowest entry representing the simplex $CA$ in $ABC$ does not come from the same Morse set. 

The other main departure from the classical persistence algorithm is in the way we output the
result which is not the entire reduced matrix but only a submatrix. A column in the reduced
matrix is called \emph{targetable} if its simplex appears as a lowest entry in a homogeneous column. For example, the columns $B$, $C$, $D$ are targetable as they appear as the lowest entry in the homogeneous columns $AB$, $BC$, and $DA$ respectively. All homogeneous and targetable 
columns and their corresponding rows are eliminated from the output. The resulting
$5\times 5$ connection matrix is shown which represents the dynamics depicted in the Conley complex shown in Figure~\ref{fig:annulus} (right). The column for $CDA$ in the connection matrix
has two $1$'s in rows $CA$ and $CD$ signifying that the dynamics flow from the repeller triangle
$CDA$ to 
the exisiting saddle edge $CA$ and the new saddle edge $CD$ arising 
from breaking the attracting periodic orbit. 
Similarly, the column $ABC$ has a single $1$ in row $CA$ 
but no $1$ in row $CD$ indicating that dynamics 
can flow from the triangle $ABC$ only to the repeller edge $CA$.

Applying the same algorithm on the boundary matrix shown in the bottom row, we arrive 
at a different connection matrix. Observe the difference in the entries for columns
$ABC$ and $CDA$ which corresponds to the difference in the Conley complexes shown in
Figure~\ref{fig:annulus}.

\section{Algebraic formulation}
\label{sec:prelim}

Although connection matrices were constructed by R.\ Franzosa in \cite{Fr1986,Fr1989}  
as a tool facilitating the detection of heteroclinic connections between invariant sets of dynamical systems,
they may be decoupled from dynamics, defined purely algebraically  and studied as a part of homological algebra. 
A formal setting is that the
input multivector field $\mv$ as a partition of a simplicial complex $K$ is
presented with a Morse decomposition. This decomposition is viewed as a 
\emph{chain complex} with a \emph{boundary map} that is \emph{filtered} with a poset $\mathbb{P}=(P,\leq)$ over which the Morse decomposition is supported. This input
chain complex is successively converted to other chain complexes, the last one
being the source of a Conley complex. Successive chain complexes are always kept connected
with isomorphic \emph{chain maps} presented as matrices w.r.t. some appropriate bases.
The Conley complex is extracted from the last chain complex by an implicit
chain homotopy that is implemented by a simple removal of a certain submatrix 
of the matrix representing the final chain map. 
Using homology groups with a field coefficient for each individual simplicial complex in the
\emph{filtered chain complexes}, we obtain \emph{graded vector spaces} connected by 
\emph{filtered linear
maps} induced by chain maps.
Our algebraic framework 
builds on these concepts/objects, which have
been used in earlier works in the context of connection matrix~\cite{HMS2021a,MW2021b}.

Before presenting this formal, purely algebraic definition of connection matrix,
we present some general assumptions, notations, and
some basic facts from homological algebra. 

\subsection{Algebraic preliminaries}
Throughout the paper we consider only finite dimensional vector spaces over the field $\ZZ_2$ and matrices with $\ZZ_2$ entries.
Let $V$ be an $n$-dimensional vector space over  $\ZZ_2$.
Consider the ordered set 
$\II_n:=(1,2,\ldots,n)$.
Since in the paper the order of vectors in a basis $B=(b_1,b_2,\ldots,b_n)$ of $V$ matters,
we consider the basis as a map $b:\II_n\ni i\mapsto b_i\in V$.
Given a basis $B$, we denote the associated scalar product by $\scalprod{\cdot,\cdot}$.
More precisely,  
$\scalprod{\cdot,\cdot}$ is
the bilinear form $V\times V\to\ZZ_2$ defined on generators $b,b'\in B$ by
\[
\scalprod{b,b'}=\begin{cases}
                 1 & \text{ if $b=b'$,}\\
                  0 & \text{ otherwise}
                \end{cases}
\]
and linearly extended to $V\times V$.
An example of an $n$-dimensional vector space over $\ZZ_2$ is the coordinate space $\ZZ_2^n$.
The {\em canonical basis} $\bE_n=\{\be_1,\be_2,\ldots,\be_n\}$ consists of vectors $\be_i$ in which all coordinates except the $i$th are zero.

Given a matrix $A$, we denote its $i$th row by $A[i,\cdot]$, its $j$th column by $A[\cdot,j]$
and its entry in $i$th row and $j$th column by $A[i,j]$.
By $\low_A(j)$ we denote the row index of the lowest $1$ in $A[\cdot,j]$ if the column is nonzero and we set $\low_A(j)=0$ otherwise. 
If $i:=\low_A(j)$ is nonzero,  we say that column  $A[\cdot,i]$ is the {\em target} of column  $A[\cdot,j]$. 


Let $V'$ be an $n'$-dimensional vector space over  $\ZZ_2$ with basis $B'=\{b'_1,b'_2,\ldots,b'_{n'}\}$.
Given a linear map $h:V\to V'$, its {\em matrix in bases} $B$ {\em and} $B'$ is the matrix $(h_{ij})$ where
\begin{equation}
\label{eq:hij}
h_{ij}:=\scalprod{hb_j,b'_i} \text{ for $i\in\II_{n'}$ and $j\in\II_n$.}
\end{equation}
Conversely, given a $\ZZ_2$-matrix $A=(a_{ij})$ with $n'$ rows and $n$ columns 
we have a linear map $h_{B,A,B'}:V\to V'$ defined on  basis $B=(b_1,b_2,\ldots,b_n)$
by $h_{B,A,B'}(b_j):=\sum_{i=1}^{n'}a_{ij}b'_i$. 
Then, clearly, $A$ is the matrix of $h_{B,A,B'}$.
If bases $B$ and $B'$ are clear from context, we adopt the simplified notation $h_A$.

Given a fixed basis $B$ of $V$, we identify $V$ with the coordinate space $\ZZ_2^n$ via isomorphism $h_{B,I_n,\bE_n}$
where $I_n$ denotes an $n\times n$ identity matrix. 
Under this identification, composition of linear maps corresponds to multiplication of matrices.

We recall (see, for instance, \cite[Section 3.1]{KaMiMr2004}) 
that for $1\leq i<j\leq n$ the operation of adding $i$th column to $j$th column followed by 
adding $j$th row to $i$th row in an $n\times n$ $\ZZ_2$-matrix $A$
results in a new matrix $A'$ such that 
\begin{equation}
\label{eq:AE}
E_{i,j} A'=A E_{i,j}
\end{equation}
where $E_{i,j}$ (illustrated below) is the sum of the identity matrix and the matrix whose all entries are zero except a one in $i$th row and $j$th column.
We call it the {\em column/row addition} matrix.
      \begin{equation*}
      \label{eq:cl-matrix}
        E_{i,j}:=\begin{array}{cc}
                     \left[
                     \begin{array}{ccccccccc}
                       1 & & & & & & & &  \\
                         &.& & & & & & &  \\
                         & &1&0&.&0&1& &  \\
                         & & &1& & &0& &  \\
                         & & & &.& &.& &  \\
                         & & & & &1&0& &  \\
                         & & & & & &1& &  \\
                         & & & & & & &.&  \\
                         & & & & & & & & 1\\
                     \end{array}
                     \right]
                     &
                     \begin{array}{c}
                          \\
                          \\
                         (i)\\
                          \\
                          \\
                          \\
                         (j)\\
                          \\
                          \\
                     \end{array}
                   \end{array}
      \end{equation*}

We recall that a {\em chain complex} with $\ZZ_2$ coefficients is a pair $(C,d)$ where $C=\bigoplus_{q\in\ZZ}C_q$
is a $\ZZ_2$ vector space with gradation $C=(C_q)_{q\in \mathbb{Z}}$ 
and $d:C\to C$ is a linear map satisfying $d(C_q)\subset C_{q-1}$ and $d^2=0$.
Given another such chain complex $(C',d')$, a {\em chain map} $\varphi:(C,d)\to(C',d')$ is a linear map such that 
$\varphi(C_q)\subset C'_q$ and $d'\varphi=\varphi d$.
A chain map is a {\em chain isomorphism} if it is an isomorphism as a linear map. 
Two chain maps $\varphi,\psi:(C,d)\to(C',d')$ are {\em chain homotopic} if there exists a linear map $S: C\to C'$
satisfying $S(C_q)\subset C'_{q+1}$ and $\psi-\varphi=d' S+S d$. Such an $S$ is called a {\em chain homotopy}.
Two chain complexes $(C,d)$, $(C',d')$ are {\em chain homotopic} if there exist chain maps $\varphi:(C,d)\to(C',d')$
and $\psi:(C',d')\to(C,d)$ such that $\psi\varphi$ is chain homotopic to $\id_C$ and $\varphi\psi$ is chain homotopic to $\id_{C'}$.
In what follows, we consider only finite dimensional chain complexes that allow
us to work with finite bases.

\subsection{Algebraic connection matrices}
Let $P$ be a finite set.
By a $P$-{\em gradation} of a finite dimensional $\ZZ_2$-vector space $V$ we mean 
the collection $\setof{V_p\mid p\in P}$ of subspaces of $V$ such that $V=\bigoplus_{p\in P}V_p$.
We call such a vector space with a given $P$-gradation a $P$-{\em graded vector space}.
We say that basis $B=(b_1,b_2,\ldots,b_n)$ is $P$-{\em graded}
if $B\subset\bigcup_{p\in P}V_p$. 
Note that for each basis vector $b_i\in B$ there is exactly one $p\in P$ such that $b_i\in V_p$. 
We denote this $p$ by $\pgrad{P}{b_i}$.
Clearly, every $P$-{\em graded vector space} admits a $P$-graded basis. 

Given a $P$-graded vector space $V=\bigoplus_{p\in P}V_p$ 
we denote by $\iota^V_{q}:V_q\to V$ and $\pi^V_{p}:V\to V_p$ respectively the inclusion and projection homomorphisms.
Given another $P$-graded vector space $V'=\bigoplus_{p\in P}V'_p$ 
we identify linear map $h:V\to V'$ with the matrix $[h_{pq}]_{p,q\in P}$ of partial linear maps
$h_{pq}:V_q\to V'_p$ where $h_{pq}:= \pi^{V'}_p\circ h\circ \iota^V_q$.

Consider now a fixed finite poset $\PP:=(P,\leq)$ and a linear map $h:V\to V'$.
We say that $h$ is $\PP$-{\em filtered} or briefly {\em filtered} when $\PP$ is clear from the context if
\begin{equation}
\label{eq:filt-homo}
h_{pq}\neq 0 \implies p\leq q.
\end{equation}
The following proposition is straightforward.
\begin{prop}
\label{prop:extension}
If $h$ is $\PP$-filtered then it is $\PP'$-filtered for any poset $\PP'=(P,\leq')$ with partial order $\leq'$ extending $\leq$.
\end{prop}

Given a fixed $P$-graded basis $B$ of $V$ and $B'$ of $V'$, we say that an  ${n'}\times n$ matrix $A$ is $\PP$-{\em filtered} w.r.t. $B$, $B'$ 
or briefly {\em filtered} when $\PP$, $B$ and $B'$ are clear from the context and
if the linear map $h_A:V\to V'$ is $\PP$-{\em filtered}.
One can easily verify the following proposition.
\begin{prop}
\label{prop:filt-matr}
The product of $\PP$-filtered matrices is $\PP$-filtered and the inverse of an invertible
$\PP$-filtered matrix is $\PP$-filtered.
\qed
\end{prop}
By a $\PP$-{\em filtered} chain complex we mean
a chain complex $(C,d)$ with field coefficients and a given gradation $C=\bigoplus_{p\in P}C_p$
such that the boundary homomorphism $d$ is $\PP$-{\em filtered}.
Given another $\PP$-filtered chain complex  $(C',d')$ we define a $\PP$-{\em filtered chain map} 
$\varphi:(C,d)\to (C',d')$ as a chain map which is also $\PP$-filtered as a homomorphism.
One can easily check that a composition of filtered chain maps is a filtered chain map. 
We say that $\varphi$ is a {\em filtered chain isomorphism} if it is a filtered chain map which is also an isomorphism. 
Note that trivially the identity homomorphism on $(C,d)$, denoted $\id_C$, is a filtered chain isomorphism.
Two filtered chain complexes $(C,d)$ and  $(C',d')$  are {\em filtered chain isomorphic} if there exist filtered chain maps
$\varphi:(C,d)\to (C',d')$ and $\varphi':(C',d')\to (C,d)$ such that $\varphi'\circ\varphi=\id_C$ 
and $\varphi\circ\varphi'=\id_{C'}$.

Two filtered chain maps $\varphi,\varphi':(C,d)\to (C',d')$ are said to be {\em filtered chain homotopic} if there exists
a chain homotopy joining  $\varphi$ with $\varphi'$ which is also filtered as a homomorphism.
Two filtered chain complexes $(C,d)$ and  $(C',d')$  are {\em filtered chain homotopic} if there exist filtered chain maps
$\varphi:(C,d)\to (C',d')$ and $\varphi':(C',d')\to (C,d)$ such that $\varphi'\circ\varphi$ is filtered chain homotopic to $\id_C$
and $\varphi\circ\varphi'$ is filtered chain homotopic to $\id_{C'}$.

Following \cite{HMS2021,MW2021b} we define {\em Conley complex} of a filtered chain complex $(C,d)$ as any filtered chain complex $(\bar{C},\bar{d})$
which is filtered chain homotopic to $(C,d)$ and satisfies $\bar{d}_{pp}=0$ for all $p\in\PP$.
The following theorem may be easily obtained as a consequence of results in~\cite[Theorem 8.1, Corollary 8.2]{RoSa1992} and~\cite[Proposition 4.27]{HMS2021}.
\begin{theorem}
\label{thm:connection-matrix}
For every finitely generated chain complex $(C,d)$ there exists a Conley complex and any two Conley complexes of $(C,d)$ are filtered chain isomorphic. 
\qed
\end{theorem}

Theorem~\ref{thm:connection-matrix} lets us define a {\em connection matrix} of $(C,d)$ 
as the matrix of homomorphisms $[\bar{d}_{pq}]_{p,q\in\PP}$ for any Conley complex  $(\bar{C},\bar{d})$ of  $(C,d)$. We note that each homomorphism in the connection matrix is represented itself as a matrix.
Moreover, this matrix is often $1\times 1$ matrix, consisting of just one number, which is the case in our setting.

\section{Algorithm}
\label{sec:algorithm}
In this section we present algorithm {\sc ConnectMat} which computes a connection matrix of a filtered chain complex $(C,d)$.
The algorithm has a reduction phase followed by a simple
extraction step whose correctnesses are argued in subsections~\ref{sec:reductions}
and~\ref{sec:extractions} respectively.
\subsection{Reductions}
\label{sec:reductions}
We need to choose a special basis of $C$, given by the following definition. 

\begin{defn}
\label{def:admissible-basis}
{\em
  We say that a $\PP$-graded basis $B$ of chain complex $(C,d)$ is $d$-{\em admissible} if there
  exists an  extension 
  of the partial order $\leq$ in $\PP$ to a linear order $\leqlin$ such that 
\begin{equation}
\label{eq:grad-monotone}
i\leq j\implies \pgrad{P}{b_i}\leqlin \pgrad{P}{b_j} \text{ for $i,j\in\II_n$.}
\end{equation}
  and the matrix $A$ of the boundary operator $d$ in this basis is strictly upper triangular, i.e., 
\begin{equation}
\label{eq:admissible}
A[i,j]\neq 0 \implies i<j.
\end{equation}
}
\end{defn}

As an immediate consequence of \eqref{eq:grad-monotone} we get that for every $d$-admissible basis $B$ 
\begin{equation}
\label{eq:filtred-reverse}
\pgrad{P}{b_i}<\pgrad{P}{b_j} \implies i<j. 
\end{equation}

\begin{prop}
\label{prop:dab-exists}
Filtered Chain complex $(C,d)$ admits a $d$-admissible basis.
\end{prop}
\proof
  A linear extension of $\leq$ always exists
  and may be constructed by algorithms known as topological sorting, for instance Kahn's algorithm \cite{Kahn1962TopologicalSO}.
  Let  $B$ be a $\PP$-graded basis of $C$. By reordering its elements we can assure that \eqref{eq:grad-monotone} is satisfied.
  Then \eqref{eq:admissible} follows from \eqref{eq:filtred-reverse} if $\pgrad{P}{b_i}\neq\pgrad{P}{b_j}$
  and when $\pgrad{P}{b_i}=\pgrad{P}{b_j}$, \eqref{eq:admissible} may be achieved by a suitable rearrengement of basis elements 
  with the same grade. 
\qed

In the sequel we assume that $B$ is a fixed $d$-admissible basis of chain complex $(C,d)$ and $\PPlin:=(P,\leqlin)$ stand for the associated linear
extension of poset $\PP$.
It follows from \eqref{eq:admissible} that
\begin{equation}
\label{eq:lowj-j}
\low_A(j)<j \text{ for all $j\in\II_n$.}
\end{equation}


\begin{algorithm}
    \caption{{\sc{ConnectMat}}}\label{alg:connectmat}
    \KwData{An $n\times n$ matrix $A$ of a filtered boundary homomorphism $d$} 
    \KwResult{A connection matrix}
    \For(){$j:=1$ to $n$}{ 
        \For{$i:=\low_{A}(j)$ \text{down to} $1$}{
            \If{$A[i,j]=1$}{ 
                $S:=\setof{s\in\II_n\mid \text{$s\neq j$ \& $\low_{A}(s)=i$ \& $A[\cdot,s]$ is homogeneous}}$\;
                \If{$S\neq\emptyset$}{ 
                    $s:=\min S$\;
                    add column $A[\cdot,s]$ to column $A[\cdot,j]$\;
                    add row $A[j,\cdot]$ to row $A[s,\cdot]$\; 
                }
            }
        }
    }
$J:=\II_n\setminus J_{h}(A)\setminus J_{t}(A)$\;
\Return{$A$ restricted to columns and rows with indices in $J$}
\end{algorithm}
The algorithm {\sc ConnectMat} computing a connection matrix of $\PP$-filtered chain complex $(C,d)$ 
is presented in Algorithm~\ref{alg:connectmat}. 
In order to explain the algorithm and discuss its features we first make some assumptions and introduce some notation. 
We assume that {\sc ConnectMat}  
takes as input a filtered matrix $A$ of the boundary homomorphism $d$ in a fixed $d$-admissible 
basis $B=(b_1,b_2,\ldots b_n)$. We denote the associated linear extension of the partial order in $\PP$ by $\leqlin$.
As we already mentioned in Section~\ref{sec:prelim}, fixing basis $B$ lets us
identify $(C,d)$ with ($\ZZ_2^n,A)$ where $n$ is the dimension of $C$ and $A$ is the matrix of the boundary homomorphism $d$
in basis $B$.  
For $i\in\II_n$ we set $\nu(i):=\pgrad{P}{b_i}$. By \eqref{eq:grad-monotone} we have
\begin{equation}
\label{eq:nu-monotone}
i\leq j\implies \nu(i)\leqlin \nu(j) \text{ for $i,j\in\II_n$.}
\end{equation}
We call a column $j$ of matrix $A$ {\em homogeneous} if it is nonzero and $\nu(j)= \nu(\low_A(j))$.
We call a column $j$ {\em targetable} if it is the target of a homogeneous column.
We denote the set of homogeneous and targetable columns of $A$ respectively by $J_h(A)$ and $J_{t}(A)$.

Let $K$ be the number of times the inner \kwfor loop is executed by the algorithm. 
Clearly, $K$ is finite.
Let $A=A_0$ be the initial matrix
before entering the $\kwfor$ loop and $A=A_k$ be the matrix after
the $k$th execution of the inner \kwfor loop for $k=1,2,\ldots,K$.
Let $i_k$, $j_k$ and $s_k$ denote the corresponding values of the variables $i$, $j$ and $s$ respectively.
Let $d_k:=h_{A_k}$ and set $E_k:=E_{s_k,j_k}$ if $s_k\not= 0$ and the identity matrix
otherwise. 
We get from \eqref{eq:AE} that for $k=1,2,\ldots, K$  
\begin{equation}
\label{eq:EAE}
    A_k=E_k^{-1}A_{k-1}E_k.
\end{equation}

\begin{prop}
\label{prop:Ak}
  For every $k\in\II_N$ matrices $A_k$ and $E_k$ are filtered. Moreover, $(C,d_k)$ is a filtered chain complex
  which is filtered chain isomorphic to $(C,d)$.
\end{prop}
\proof
Clearly, $A_0^2=0$, because $A_0$ equals the matrix on input which is the matrix of a boundary operator.
Therefore, by induction argument we get $A_k^2= E_k^{-1}A_{k-1}^2 E_k=0$ which shows that $d_k=h_{A_k}$ is a boundary operator.
Moreover, we get from \eqref{eq:EAE} 
that  $e_k:=h_{E_k}$ is a chain map
and by definition it is an isomorphism from  $(C,d_{k-1})$ to $(C,d_k)$.
Therefore, it is a chain isomorphism.

We will prove by induction on $k$ that $A_k$ is filtered. 
Again, $A_0$ as the matrix on input is filtered. 
Assume $A_{k-1}$ is filtered for some $k$ satisfying $1\leq k \leq K$. We first show that $E_k$ is filtered.
To this end assume that $(h_{E_k})_{pq}\neq 0$ for some $p,q\in\PP$. In view of the structure of matrix $E_k$ either $(h_{E_k})_{pq}$ contains a diagonal entry
which immediately implies $p=q$, or it contains the unique off-diagonal nonzero entry in $s_k$th row and $j_k$th column. 
Consider the latter case. Then $p=\nu(s_k)$ and $\nu(j_k)=q$. 
Since variable $A$ is modified,
the condition in the \kwif statement is satisfied and the properties of the set
$S$ hold. 
Hence, we have $A_{k-1}[i_k,j_k]=1$ and $\low_{A_{k-1}}(s_{k})=i$. We also know that column $s_k$ is homogeneous. 
Since $A_{k-1}$ is filtered by induction assumption, we get $\nu(i_k)\leq \nu(j_k)$.
Combining all these properties we get 
\[
p=\nu(s_k)=\nu(\low_{A_{k-1}}(s_{k}))=\nu(i_k)\leq\nu(j_k)=q.
\]
Hence, in both bases we get $p\leq q$, which proves that $E_k$ is filtered
and in view of Proposition~\ref{prop:filt-matr} we get from \eqref{eq:EAE} that $A_k$  is filtered.
Thus, we have filtered chain complexes $(C,d_k)$ for $k=0,1,\ldots, K$ and all these chain complexes are filtered chain isomorphic.
\qed

\begin{prop}
\label{prop:low}
  For every $j\in\II_n$ the $j$th column of matrix $A$ does not change after completing the $j$th pass of the outer \kwfor loop.
\end{prop}
\proof
  Fix a $j\in\II_n$. 
  Note that the only modifications of matrix $A$ happen in the addition of columns and rows inside the inner $\kwfor$ loop.
  Consider column and row additions on $k$th pass of the outer \kwfor loop with $k>j$.
  Clearly, column additions modify $k$th column but not $j$th column.
  A row addition  on $k$th pass of the outer \kwfor loop adds $k$th row to another row.
  Note that by \eqref{eq:lowj-j} we have $\low_A(j)<j<k$. Therefore, $A[k,j]=0$, which means that such an addition of rows cannot
  modify $j$th column as well. 
\qed

\begin{prop}
\label{prop:Aij}
The value of $A[i,j]$ settled after completing the 
 inner \kwfor loop with the index variables {\em $i$} and {\em $j$} remains unchanged till the end of the algorithm.
\end{prop}
\proof
  Assume to the contrary that such a change happens on executing the inner \kwfor loop when  $j=j'$ and $i=i'$.
  By Proposition~\ref{prop:low} we must have $j'=j$ and since the inner \kwfor loop which is 
  a downwards loop, we have $i'<i$.
  The change of $A[i,j]$ may happen only via addition of columns or rows inside \kwif statement.
  Therefore, for some $s'$, we have $\low_A(s')=i'<i$, which implies $A[i,s']=0$. In consequence, addition of column $A[\cdot,s']$ to $A[\cdot,j]$ cannot change $A[i,j]$. 
  Addition of row $A[j,\cdot]$ to $A[s',\cdot]$ can change $A[i,j]$ only if $s'=i$ and $A[j,j]=1$ which clearly is not the case. This contradiction proves our claim. 
\qed

\begin{prop}
\label{prop:reducible}
After completing
the $j$th pass (for any $j\in \II_n$) of the outer for loop, if
$A[i,j] = 1$ for some $i\leq \low_A(j)$, then there is no 
  homogeneous column $A[\cdot,s]$ with $s\neq j$ and $\low_A(s)=i$.
\end{prop}
\proof 
  Assume on the contrary that for some $j\in\II_n$ after completing $j$th pass of the outer \kwfor loop there is an $i\leq \low_A(j)$ such that $A[i,j]=1$ 
  and a homogeneous column $A[\cdot,s]$ with $s\neq j$ and $\low_A(s)=i$.
  Consider the value of $A[i,j]$ immediately before entering the inner \kwfor loop with index variables $i$
  and $j$ for the inner and outer $\kwfor$ loops respectively.
  We claim  that regardless of the value of $A[i,j]$,  we have $A[i,j]=0$
  after completing this pass of the inner \kwfor loop.
  Indeed, if the initial value is zero, then the condition of \kwif statement is not satisfied. 
  Consequently, $A[i,j]$ is not modified and remains zero after completing the inner \kwfor loop.
  If the value is one, then \kwif statement is executed and addition of columns modifies $A[i,j]$ to zero.
  This proves our claim. Now, by Proposition~\ref{prop:Aij} the value of $A[i,j]$ is also zero when the $j$th pass of the outer \kwfor loop is completed, a contradiction. 
  \qed

We say that a filtered matrix $A$ is {\em reduced} if it satisfies the following three conditions:
\begin{itemize}
   \item (R1) the map
$
  \alpha: J_h(A)\ni j\mapsto \low_{A}(j)\in J_t(A)
$
is a well defined bijection, 
   \item (R2) $J_h(A)\cap J_t(A)=\emptyset$, and
   \item (R3) if $A[i,j]=1$ for some $i\leq \low_A(j)$, then there is no $s\in J_h(A)$ with $s\neq j$ and $\low_A(s)=i$.
\end{itemize}

\begin{prop}
\label{prop:final-A-reduced}
  Matrix $A_K$, that is the final state of the variable $A$, is reduced.
\end{prop}
\proof 
To simplify the notation, denote by $A:=A_K$ the final state of variable $A$.
By definition of  $J_t(A)$, the map $\alpha$ is well defined and is a surjection. 
The fact that it is an injection follows immediately from Proposition~\ref{prop:reducible}.
This proves property (R1).

To prove (R2) assume on the contrary that there is an $i\in J_h(A)\cap J_t(A)$.
Then, there is a $j\in J_h(A)$ such that $\low_{A}(j)=i$. Since column $j$ is homogeneous, we have $\nu(i)=\nu(\low_{A}(j))=\nu(j)$.
Let $t:=\low_{A}(i)$. Since $i\in J_h(A)$, we have $t>0$. Proof of Proposition~\ref{prop:Ak} implies that $A^2=0$. Thus, we have
\[
0=A^2[t,j]=\sum_{s=1}^n A[t,s] A[s,j]=\sum_{s=t+1}^i A[t,s] A[s,j]=\sum_{s=t+1}^{i-1} A[t,s] A[s,j]+1,
\]
because $A[t,s]=0$ for $s\leq t$, $A[s,j]=0$ for $s>i$ and $A[t,i]=1=A[i,j]$. 
Hence, there is an $s_0$ such that $t<s_0<i$ and $A[t,s_0]=1=A[s_0,j]$.
In particular, since $A[t,s_0]=1$, we get from Proposiiton~\ref{prop:reducible} that there is no homogeneous column $A[\cdot,u]$ such that
$\low_A(u)=t$. However, $A[\cdot,i]$ is exactly such a column, a contradiction. 
Finally, property (R3) follows immediately from Proposition~\ref{prop:reducible}.
\qed

We say that a filtered matrix $A$ is {\em cropped} if  $J_h(A)=\emptyset$. Then, necessarily also $J_{t}(A)=\emptyset$.

\begin{prop}
\label{prop:d-cropped}
  Consider matrix $A$ of the boundary operator $d$ in a fixed $d$-admissible basis $B$ of filtered chain complex $(C,d)$.
  If $A$ is cropped, then $(C,d)$ is its own Conley complex and $A$ is a connection matrix of $(C,d)$.
\end{prop}
\proof
  Clearly $(C,d)$ is filtered chain homotopic to itself. Hence, we only need to check that $d_{pp}=0$ for $p\in\PP$.
  Assume to the contrary that there is a $p\in\PP$ such that $d_{pp}\neq 0$.
  We have $d=h_A$, therefore there are $i,j\in \II_n$ such that $\nu(i)=\nu(j)=p$ and $A[i,j]=1$.
  It follows that the $j$th column of $A$ is nonzero. Since $J_h(A)=\emptyset$, we have $\nu(\low_A(j))\neq\nu(j)$.  
  We also have $A[\low_A(j),j]=1$ which gives $\nu(\low_A(j))\leq\nu(j)$.  
  Therefore, $\nu(\low_A(j))<\nu(j)=\nu(i)$. We get from \eqref{eq:filtred-reverse} that $\low_A(j)<i$.
  But, from the definition of $\low_A$ we have $i\leq\low_A(j)$, a contradiction proving that $d_{pp}=0$ for $p\in\PP$.
\qed

\subsection{Extracting the connection matrix}
\label{sec:extractions}
Observe that the algorithm outputs the cropped matrix $A$, which by Proposition~\ref{prop:d-cropped} is a connection matrix of the 
chain complex $(\bar{C},\bar{d})$ where $A$ represents the boundary operator $\bar{d}$.
Now we show that the $(\bar{C},\bar{d})$ is filtered chain homotopic to the input
chain complex $(C,d)$. This establishes that the output matrix $A$ is the
connection matrix of $(C,d)$.

We say that a basis $B$ of a chain complex $(C,d)$ is {\em $\ZZ$-graded} if $B\subset \bigcup_{q\in\ZZ}C_q$.
Let $B=(b_1,b_2,\ldots,b_n)$ be a fixed $\ZZ$-graded basis of $C$. 
We say that a pair $(i_0,j_0)$ of indexes in $\II_n$ is a {\em reduction pair} if $b_{i_0},b_{j_0}\in B$ are such that $\scalprod{d b_{j_0},b_{i_0}}=1$.
Given a fixed reduction pair $(i_0,j_0)$ denote by $\bar{C}$ the subspace of $C$ spanned by $\bar{B}:=B\setminus\{b_{i_0},b_{j_0}\}$ 
and consider maps 
\begin{eqnarray}
\label{eq:bard}
  &&\bar{d}:\bar{C}\ni c\mapsto  d c + \scalprod{d c,b_{i_0}}d b_{j_0} + \scalprod{d c,b_{j_0}}b_{j_0}\in \bar{C}\\
\label{eq:pi}
  &&\pi:{C}\ni c\mapsto  c + \scalprod{c,b_{i_0}}d b_{j_0} + \scalprod{c,b_{j_0}}b_{j_0}\in \bar{C}\\
\label{eq:iota}
  &&\iota:\bar{C}\ni c\mapsto c + \scalprod{d c,b_{i_0}}b_{j_0}\in {C}\\
\label{eq:gamma}
  &&\gamma:C\ni c\mapsto \scalprod{c,b_{i_0}} b_{j_0}\in C.
\end{eqnarray}

The following theorem is based on results in \cite{KMS1998}.
\begin{theorem}
\label{thm:C-bar} (see \cite{KMS1998})
  The pair $(\bar{C},\bar{d})$ is a finite dimensional chain complex with coefficients in $\ZZ_2$. 
Moreover, $\iota$ and $\pi$ are chain maps and $\gamma$ is a chain homotopy such that
\[
  \iota\pi=\id_C + d\gamma + \gamma d \quad \text{ and } \quad
  \pi\iota =\id_{\bar{C}}.
\]
In particular, chain complexes $(C,d)$ and $(\bar{C},\bar{d})$ are chain homotopic.
We refer to complex $(\bar{C},\bar{d})$ as the $(i_0,j_0)$-{\em reduction} of complex $(C,d)$.
\end{theorem}
\proof
  The first part of the theorem is a special case of \cite[Theorem 1]{KMS1998}.
  The second part is implicitly contained in the proof of \cite[Theorem 2]{KMS1998}.
\qed

\begin{prop}
\label{prop:matrices}
Consider a filtered chain complex $(C,d)$ so that the matrix $(d_{ij})$ of $d$ in a $d$-admissible basis $B=(b_1,b_2,\ldots,b_n)$
is reduced.  Let $j_0\in J_h(A)$ and let $i_0:=\low(j_0)$.
Then $(i_0,j_0)$ is a reduction pair. Moreover, if $(\bar{C},\bar{d})$ is the 
$(i_0,j_0)$-reduction of complex $(C,d)$ as in Theorem~\ref{thm:C-bar}, then
the  matrices of maps given by formulas (\ref{eq:bard}-\ref{eq:gamma}) may be characterized as follows. 
The matrix $(\bar{d}_{ij})$ of $\bar{d}$ in base $\bar{B}$ satisfies for $i,j\in\bar{B}$ 
\begin{equation}
  \label{eq:bard-m}
  \bar{d}_{ij}=d_{ij}.
\end{equation}
The matrix $(\pi_{ij})$ of $\pi$ in bases $B$ and $\bar{B}$  satisfies for $i\in\bar{B}$ and $j\in B$
\begin{equation}
  \label{eq:pi-m}
  \pi_{ij}=\begin{cases}
                   \delta_{ij} & \text{if $j\neq i_0$,}\\
                   d_{ij_0} &  otherwise,
                 \end{cases}
\end{equation}
where $\delta_{ij}$ is zero unless $i=j$ when it is one.                 
The matrix $(\iota_{ij})$ of $\iota$ in bases $\bar{B}$ and  $B$  satisfies for $i\in B$ and $j\in \bar{B}$
\begin{equation}
  \label{eq:iota-m}
  \iota_{ij}=\begin{cases}
                   \delta_{ij} & \text{if $i\neq j_0$,}\\
                   d_{i_0j} & \text{otherwise.}
                 \end{cases}
\end{equation}
The matrix $(\gamma_{ij})$ of $\gamma$ in base $B$  satisfies for $i,j\in B$ 
\begin{equation}
  \label{eq:gamma-m}
  \gamma_{ij}=\begin{cases}
                   1 & \text{if $i=j_0$ and $j=i_0$,}\\
                   0 & \text{otherwise.}
                 \end{cases}
\end{equation}
\end{prop}
\proof
Since $i_0:=\low(j_0)$, clearly $(i_0,j_0)$ is a reduction pair.
To prove \eqref{eq:bard-m} observe that by \eqref{eq:hij} and \eqref{eq:bard} we have
\begin{equation*}
\bar{d}_{ij}=\scalprod{db_j,b_i}+\scalprod{db_j,b_{i_0}}\scalprod{db_{j_0},b_{i}}+\scalprod{db_j,b_{j_0}}\scalprod{b_{j_0},b_{i}}
            =d_{ij}+\scalprod{db_j,b_{i_0}}\scalprod{db_{j_0},b_{i}},
\end{equation*}
because $i\in\bar{B}$ implies $\scalprod{b_{j_0},b_{i}}=0$.
Since $j_0\in J_h(A)$, we cannot have $d_{i_0j}=\scalprod{db_j,b_{i_0}}\neq 0$, because the matrix $d_{ij}$ of $d$ is reduced.
Hence $\scalprod{db_j,b_{i_0}}\neq 0$ and \eqref{eq:bard-m} follows.

To prove \eqref{eq:pi-m} observe that 
\begin{equation*}
\pi_{ij}=\scalprod{b_j,b_i}+\scalprod{b_j,b_{i_0}}\scalprod{db_{j_0},b_{i}}+\scalprod{b_j,b_{j_0}}\scalprod{b_{j_0},b_{i}}
            =\delta_{ij}+\scalprod{b_j,b_{i_0}}\scalprod{db_{j_0},b_{i}},
\end{equation*}
because $i\in\bar{B}$ implies $\scalprod{b_{j_0},b_{i}}=0$.
Hence, if $j\neq i_0$, we get  $\scalprod{b_j,b_{i_0}}=0$ and $\pi_{ij}=\delta_{ij}$.
Since $i\neq i_0$, in the case $j = i_0$ we get $\delta_{ij}=0$ and $\pi_{ij}=\scalprod{db_{j_0},b_{i}}=d_{ij_0}$.
This proves \eqref{eq:pi-m}.

To prove \eqref{eq:iota-m} observe that 
\begin{equation*}
\iota_{ij}=\scalprod{b_j,b_i}+\scalprod{d b_j,b_{i_0}}\scalprod{b_{j_0},b_{i}}.
\end{equation*}
Hence, if $i\neq j_0$, we get $\scalprod{b_{j_0},b_{i}}=0$ and $\iota_{ij}=\scalprod{b_j,b_i}=\delta_{ij}$.
Since $j\neq i_0$, in the case $i = j_0$ we get $\delta_{ij}=0$ and $\iota_{ij}=\scalprod{d b_j,b_{i_0}}=d_{i_0j}$.
This proves \eqref{eq:iota-m}.

Finally, to prove \eqref{eq:gamma-m} observe that 
\begin{equation*}
\gamma_{ij}=\scalprod{b_j,b_{i_0}}\scalprod{b_{j_0},b_{i}}.
\end{equation*}
Therefore, $\gamma_{ij}=0$ unless $i=j_0$ and $j=i_0$ when $\gamma_{ij}=1$.
This proves \eqref{eq:gamma-m}.
\qed  \\                   
\begin{prop}
\label{prop:elem-red}
  Under assumptions of Proposition~\ref{prop:matrices}, 
the associated $(i_0,j_0)$-reduction of $(C,d)$, denoted $(\bar{C},\bar{d})$,
  is a filtered chain complex, filtered chain homotopic to $(C,d)$. Moreover, the matrix of $\bar{d}$ in basis $B\setminus\{b_i,b_j\}$
  is precisely the matrix of $d$ with $i_0$th and $j_0$th columns and rows removed.
  In particular, it is also reduced.
\end{prop}
\proof
   We get from Theorem~\ref{thm:C-bar} that $(\bar{C},\bar{d})$ is a chain complex which is chain homotopic
   to $(C,d)$. To see that $(\bar{C},\bar{d})$ is a filtered chain complex, we need to prove that $\bar{d}$ is filtered.
   For this end assume that $\bar{d}_{pq}\neq 0$ for some $p,q\in P$. Then, there are $i,j\in\II_n$ such that $b_i,b_j\in\bar{B}$,
   $\nu(i)=p$, $\nu(j)=q$, and $\bar{d}_{ij}\neq 0$. 
   Then, by \eqref{eq:bard-m}, we have $d_{ij}=\bar{d}_{ij}\neq 0$ and, since $d$ is filtered, we get $p=\nu(i)\leq\nu(j)=q$, which 
   proves that $\bar{d}$ is filtered.
   
   To prove that $(\bar{C},\bar{d})$ is filtered chain homotopic to $(C,d)$ we only needs to verify that also  $\pi$, $\iota$ and $\gamma$ 
   are filtered. This follows immediately from formulas \eqref{eq:pi-m}, \eqref{eq:iota-m}, \eqref{eq:gamma-m} and the fact that $d$ is filtered.
      Finally, formula \eqref{eq:bard-m} implies that the matrix of $\bar{d}$ coincides with that of $d$ with $i_0$th and $j_0$th 
   columns and rows removed. This immediately implies that the matrix of $\bar{d}$ is reduced.
\qed

\begin{theorem}
Given an $n\times n$ matrix of a filtered chain complex $(C,d)$ in a fixed $d$-admissible basis on input, {\sc ConnectMat}
outputs a connection matrix of $(C,d)$ in $O(n^3)$ time. 
\label{thm:main}
\end{theorem}
\proof
Let $A$ be the final matrix, that is, $A:=A_K$ and let $\bar{A}$ be the matrix on output of algorithm {\sc ConnectMat}.
First observe that by Proposition~\ref{prop:Ak} chain complex $(C,h_A)$ is filtered chain isomorphic to $(C,d)$. 
By Proposition~\ref{prop:final-A-reduced} matrix $A$ is reduced. 
We proceed by induction in the cardinality of $J_h(A)$. If $J_h(A)=\emptyset$, the conclusion follows from Proposition~\ref{prop:d-cropped}.
If the cardinality of $J_h(A)$ is positive, we first observe that by the definition of reduced matrix and Proposition~\ref{prop:reducible} 
the pair $(\low_A(j),j)$ for $j\in J_h(A)$ is an elementary reduction pair. Therefore, by the inductive assumption 
of Proposition~\ref{prop:elem-red}, {\sc ConnectMat} outputs a connection matrix. The inner $\kwfor$ loop runs in $O(n^2)$ time (similar to the persistence algorithm)
for all column and row additions incurring an $O(n^3)$
total time complexity.
\qed
\section{Conclusions}
We have presented a simple, direct, single-pass, persistence-like algorithm for computing
connection matrices from a given boundary matrix filtered according to a Morse decomposition of 
a multivector field.
Since connection matrices need not be unique in general,
it would be interesting to know whether the presented algorithm can
be made to produce all possible connection matrices. 
There is a choice in extending the partial order in $\PP$ 
to the linear order of columns in the filtered matrix, which actually may lead to different connection matrices; Figure~\ref{fig:matrixalgo} presents such an example.
It is proven in \cite{MW2021b} that, for a Forman gradient field, there is
a unique connection matrix. This seems to be related to the fact that 
each critical Morse set in such a field consists of a single simplex (single column) thus leaving no choice for reductions. 


Now that we have
a persistence-like algorithm for connection matrices, can we extend it for
    persistence of connection matrices across fields? We plan to address these
questions in future work.
\section*{Acknowledgments.} This work is partially supported by NSF grants CCF-2049010, DMS-2301360, Polish National Science Center under Opus Grant 2019/35/B/ST1/00874 
and Preludium Grant 2018/29/N/ST1/00449.
M.L. acknowledge support by the Dioscuri program initiated by the
Max Planck Society, jointly managed with the National Science Centre (Poland), 
and mutually funded by the Polish Ministry of Science and Higher 
Education and the German Federal Ministry of Education and Research.
\bibliography{refs}



\end{document}